\documentclass[12pt]{amsart}

\usepackage{latexsym, amsmath, amscd, amssymb, amsthm} 
\usepackage[T2A]{fontenc}
 \usepackage[cp1251]{inputenc}
\usepackage[ukrainian, russian, english]{babel}
\newtheorem{te}{Theorem}

 \newtheorem{lm}{Lema}

\begin{document}

\noindent

 \title{ The Poincar\'e series of the covariants of binary forms }

\author{L. Bedratyuk}\address{Khmelnitskiy national university, Insituts'ka, 11,  Khmelnitskiy, 29016, Ukraine}

\begin{abstract}
A  formula for computation of the Poincar\'e series $P_d(z)$ of the algebra of the covariants of binary $d$-form is found. 
By using it, we have computed the $P_d(z)$ for  $d \leq 20.$
\end{abstract}
\maketitle

\noindent
{\bf 1.}  Consider the natural action of the complex matrix group   $SL_2$ on  the complex vector space of binary  $d$-forms.   Let   $ \mathit{C_d}$ be corresponding algebra of    covariants. The algebra  $C_{d}$ is a finite  generated graded algebra: 

$$
C_{d}=(C_{d})_0+(C_{d})_1+\cdots+(C_{d})_i+ \cdots,
$$
and the vector spaces  $(C_{d})_i$ are all  finite dimensional.  The formal power series 
 $P_{d}(z) \in \mathbb{Z}[[z]],$
$$
P_{d}(z)=\sum_{i=0}^{\infty }\dim((C_{d})_i) z^i,
$$ 
is called the Poincar\'e series   of the algebra of   covariants  $C_{d}.$ The finite generation of the algebra of   covariants implies  that its Poincar\'e series is an expansion  of some  rational function.  We consider here the problem of
computing efficiently this rational function.

  The  Poincar\'e series  calculations were an  important  object of research  in  classical invariant theory of the 19th century.
For  the cases  $d\leq 10,$ $d=12$ the  series  $P_{d}(z)$  were calculated by Sylvester and  Franklin, see  \cite{SF}, \cite{Sylv-12}.
Relatively recently, Springer \cite{SP} set the explicit formula for computing the Poincar\'e  series of the algebras of invariants of the binary $d$-forms. This formula has been used by Brouwer and Cohen  \cite{BC} for  $d\leq 17$ and also by Littelmann and Procesi  \cite{LP} for even  $d\leq 36.$

 In the paper we  have  found a Sylvester-Cayley-type formula  for calculating of $\dim (C_{d})_i$   and Springer-type formula  for calculation of $P_d(z).$ By using the  formula, the series $P_d(z)$  is calculated  for $d\leq 30.$

{\bf 2.} To begin with, we give a proof of the Sylvester-Cayley-type  formula for covariants of   binary form. 

Let  $V\cong \mathbb{C}^2$ be standard  two-dimensional representation of Lie algebra  $\mathfrak{sl_{2}}.$ The irreducible representation   $V_d=\langle v_0,v_1,...,v_d \rangle,$ $\dim V_d=d+1$ of the  algebra $\mathfrak{sl_{2}}$ is the symmetric $d$-power of the standard representation  $V=V_1,$  i.e. $V_d=S^d(V),$  $V_0 \cong  \mathbb{C}.$  
The basis elements    $ \left( \begin{array}{ll}  0\, 1 \\ 0\,0 \end{array} \right),$ $ \left( \begin{array}{ll}  0\, 0 \\ 1\,0 \end{array} \right)$, $ \left( \begin{array}{ll}  1 &  \phantom{-}0 \\  0 &-1 \end{array} \right)$ of the algebra    $\mathfrak{sl_{2}}$ act on    $V_d$  by the derivations  $D_1, D_2, E$ : 
$$
D_1(v_i)=i\, v_{i-1},  D_2(v_i)=(d-i)\,v_{i+1}, E(v_i)=(d-2\,i)\,v_i.
$$

The action of   $\mathfrak{sl_{2}}$  is extended to action on the symmetrical algebra  $S(V_d)$ in the natural way. 

Let   $\mathfrak{u}_{2}$ be  the maximal unipotent subalgebra of $\mathfrak{sl}_{2}.$ The algebra  $\mathcal{S}_{d}$
$$
\mathcal{S}_{d}:= \displaystyle{ S(V_d)^{\mathfrak{u_{2}}}}=\{ v \in S(V_d)|  D_1(v)=0 \},
$$
is called the algebra of {\it seminvariants}  of the binary form  of degree $d.$ For any element $v \in \mathcal{S}_{d}$ a natural number $s$ is called the {\it order} of the element $v$ if the number $s$ is the smallest natural number such that \begin{equation*}D_2^s(v) \ne 0, D_2^{s+1}(v) = 0.\end{equation*}
It is  clear that any seminvariant  $v \in S_d$  of  order $i$ is the highest weight vector  for an  irreducible $\mathfrak{sl_{2}}$-module   of the dimension $i+1$ in $S(V_d).$

The classical Robert's theorem \cite{Rob} implies  an isomorphism of the algebras covariants and seminvariants. Thus, it is  enough to compute the Poincar\'e  series of the algebra $S_d.$ 

The algebra  $S(V_d)$  is graded 
$$
S(V_d)=S^0(V_d)+S^1(V_d)+\cdots +S^n(V_d)+\cdots,
$$
and each   $S^n(V_d)$ is the complete reducibly
 representation of the Lie algebra  $\mathfrak{sl_{2}}.$
Thus, the following decomposition  holds
$$
S^n(V_d) \cong \gamma_d(n,0) V_0+\gamma_d(n,1) V_1+ \cdots +\gamma_d(n,d\,n) V_{d\,n},  \eqno{(*)}
$$
here  $\gamma_d(n,k)$ is  the  multiplicity of the representation  $V_k$  in the decomposition of  $S^n(V_d).$ On the other hand, the multiplicity  $\gamma_d(n,i)$  of the  representation  $V_i$ is  equal to the number of linearly independent homogeneous seminvariants of  degree $n$   and order $i$  for the  binary $d$-form.
This argument proves 
\begin{lm}
$$
\dim (C_d)_n=\gamma_d(n,0)+\gamma_d(n,1) + \cdots +\gamma_d(n,d\,n).
$$
\end{lm}

The set of weights ( eigenvalues of the operator $E$) of a representation  $W$ denote by  $\Lambda_{W},$  in particular, $\Lambda_{V_d}=\{-d, -d+2, \ldots, d \}.$ 

 A  formal sum 
$$
{\rm Char}(W)=\sum_{k \in \Lambda_{W}} n_W(k) q^k,
$$
is called the character   of a representation  $W,$  
here   $n_W(k)$ denotes the   multiplicity  of the weight $k \in \Lambda_{W}.$
Since, a   multiplicity of any weight of the irreducible representation $V_d$  is  equal to 1, we have  
$$
{\rm Char}(V_d)=q^{-d}+q^{-d+2}+\cdots+q^{d}.
$$

The  character   $ {\rm Char}(S^n(V_d))$ of the representation  $S^n(V_d)$  equals   $H_n(q^{-d},q^{-d+2},\ldots,q^{d}),$ see \cite{FH},  where   $H_n(x_0,x_1,\ldots,x_d)$ is  the complete symmetrical function    $$H_n(x_0,x_1,\ldots,x_d)=\sum_{|\alpha|=n} x_0^{\alpha_0} x_1^{\alpha_1} \ldots x_d^{\alpha_d} , |\alpha|=\sum_i \alpha_i.$$

By replacing    $x_k$ with $q^{d-2\,k},$ $k=0,\ldots, d,$  we  obtain the specialized expression for the character of   ${\rm Char}(S^n(V_d)):$ 
$$
{\rm Char}(S^n(V_d))= \sum_{|\alpha|=n} (q^d)^{\alpha_0} (q^{d-2\cdot 1})^{\alpha_1} \ldots (q^{d-2\,d})^{\alpha_d} =
$$
$$
= \sum_{|\alpha|=n} q^{d\,n-2 (\alpha_1+2\alpha_2+\cdots + d\, \alpha_d)}=\sum_{k=0}^{d\,n} \omega_d(n,k) q^{d\,n-2\,k},
$$
here   $\omega_d(n,k) $  is the number nonnegative integer solutions of the equation $$\alpha_1+2\alpha_2+\cdots + d\, \alpha_d=\displaystyle \frac{d\,n-k}{2}$$  on the assumption that  $ \alpha_0+\alpha_1+\cdots + d\, \alpha_d=n.$ In particular, the coefficient of $q^0$ (the     multiplicity  of zero weight ) is  equal to  $ \omega_d(n,\frac{d\,n}{2}),$ and the coefficient of  $q^1$ is  equal  $ \omega_d(n,\frac{d\,n-1}{2}).$

On the other hand, the decomposition $(*)$ implies the  equality for the characters:  
$$
{\rm Char}(S^n(V_d))=\gamma_d(n,0) {\rm Char}(V_0) +\gamma_d(n,1) {\rm Char}(V_1)+ \cdots +\gamma_d(n, d\,n) {\rm Char}(V_{d\,n}).
$$
We can summarize what we have shown so far in 
\begin{te} $$\dim (C_d)_n= \omega_d \left(n,\frac{d\,n}{2}\right)+ \omega_d \left(n,\frac{d\,n-1}{2}\right)
.$$
\end{te}
\begin{proof}
The zero weight appears  once  in any  representation $V_k,$  for even $k$,  therefore 
$$
 \omega_d \left(n,\frac{d\,n}{2}\right)=\gamma_d(n,0)  +\gamma_d(n,2) + \cdots +\gamma_d(n, 4)+...
$$
The weight   $1$ appears  once in any  representation   $V_k,$  for odd $k$,  therefore 
$$
 \omega_d \left(n,\frac{d\,n-1}{2}\right)=\gamma_d(n,1)  +\gamma_d(n,3) + \cdots +\gamma_d(n, 5)+...
$$
Thus 
$$
\omega_d \left(n,\frac{d\,n}{2}\right)+ \omega_d \left(n,\frac{d\,n-1}{2}\right)=\gamma_d(n,0)  +\gamma_d(n,1) + \cdots +\gamma_d(n, d\,n).
$$
\end{proof}

\noindent
{\bf 3.} It well-known  that  the number
 $\omega_d \left(n,\frac{d\,n}{2}\right)$  of   non-negative integer solutions of the following system
$$
\left \{
\begin{array}{c}
\displaystyle \alpha_1+2\alpha_2+\cdots + d\, \alpha_d=\displaystyle \frac{d\,n}{2},\\
\\
\displaystyle \alpha_0+\alpha_1+\cdots + \alpha_d=n
\end{array}
\right.
$$
is equal to the coefficient of $\displaystyle t^n z^{\frac{d\,n}{2}} $  of the  expansion   of the  series

$$
f_{d}(t,z)=\frac{1}{(1-t )(1-t\,z)\ldots (1-t\,z^d)}.
$$
Denote it in such a way:  $\omega_d \left(n,\frac{d\,n}{2}\right)=\left[  t^n (z)^{\frac{d\,n}{2}}\right](f_{d}(t,z)).$  It is  clear that 
$$\omega_d \left(n,\frac{d\,n}{2}\right)=\left[  t^n z^{d\,n}\right](f_{d}(t,z^2))=\left[ (t\,z^d)^n\right](f_{d}(t,z^2)).$$

Similarly, the number  $\omega_d \left(n,\frac{d\,n-1}{2}\right)$ of   non-negative integer solutions of the following system
 
$$
\left \{
\begin{array}{c}
\displaystyle \alpha_1+2\alpha_2+\cdots + d\, \alpha_d=\displaystyle \frac{d\,n-1}{2},\\
\\
\displaystyle \alpha_0+\alpha_1+\cdots + \alpha_d=n
\end{array}
\right.
$$
equals  $$\left[t^n (z)^{\frac{d\,n-1}{2}}\right](f_{d}(t,z))=\left[ t^n z^{d\,n-1}\right]( f_{d}(t,z^2))=\left[ (t z^d)^n\right]( z\,f_{d}(t,z^2)).$$
  
Thus,  the  following statement holds
\begin{te}  The number  $\dim (C_d)_n$ of   linearly independet seminvariants (and  covariants) of degree $n$ for the  binary  $d$-form is calculated by the formula
$$
\dim (C_d)_n=\left[ (t z^d)^n\right]\left( \frac{1+z}{(1-t )(1-t\,z^2)\ldots (1-t\,z^{2\,d})}\right).
$$
\end{te}

{\bf 4.} Let us prove  Springer-type  formula for the Poincar\'e  series  $P_d(z)$ of the  algebra covariants of  the binary $d$-form.  Consider the $\mathbb{C}$-algebra $\mathbb{C}[[t,z]]$   of  formal   power series.
For arbitrary   $m,n \in \mathbb{Z^+}$  define  $\mathbb{C}$-linear function
$$ \Psi_{m,n}:\mathbb{C}[[t,z]] \to \mathbb{C}[[z]],$$ $ m,n  \in \mathbb{Z}^{+} $ in the  following  way
$$
\Psi_{n_1,n_2}\bigl(t^{m_1} z^{m_2}\bigr)=\left \{ \begin{array}{l} z^{s}, \text{  if  } \displaystyle  \frac{m_1}{n_1}=\frac{m_2}{n_2}= s \in \mathbb{N}, \\ 1,  \text{  if  } m_1=m_2=0, \\ 0, \text{ otherwise. } \end{array} \right. 
$$
Then  for arbitrary series  
$$
A=a_{0,0}+a_{1,0} t +a_{0,1}z+ a_{2,0}t^2+\cdots,
$$
we  get
$$
\Psi_{n_1,n_2}(A)=a_{0,0}+a_{n_1,n_2} z+a_{2\,n_1,2\,n_2}z^2+\cdots  .
$$

Define by $\varphi_n$  the restriction of 
$\Psi_{m,n}$ to  $\mathbb{C}[[z]],$ namely 

$$
\varphi_{n}\bigl(z^m\bigr):=\left \{ \begin{array}{l} z^{\frac{m}{n}}, \text{  if  } m= 0 \pmod n,  \\ 0, \text{ if } m \neq 0 \pmod n,\\ 1,   \text {    for }   m=0. \end{array} \right. 
$$

It is clear that for arbitrary series  $$A=a_0+a_1 z+a_2 z^2+ \cdots ,$$ we  obtain
$$
\varphi_{n}(A)=a_0+a_n z+a_{2\,n} z^2+\cdots+a_{s\, n} z^s+\cdots .
$$

In  some  cases   calculation of the functions $ \Psi$  can  be reduced   to  calculation of  the functions    $\varphi$.  The following  statements hold: 
 
\begin{lm}
$$
\begin{array}{ll}
 (i) &  \text{ For  } h(t,z) \in \mathbb{C}[[t,z]]) \text{  we  have  } \displaystyle \Psi_{1,n}(h(t,z)) =\frac{1}{2 \pi i} \oint_{|z|=1} h\left(\frac{t}{z^n},z\right) {\frac{dz}{z} \Bigl |_{t=z}};
\\
(ii) & \text{ for }  R(z) \in \mathbb{C}[[z]]   \text{    and for }  m, n, k  \in \mathbb{N}     \text{  holds} \\

&  \displaystyle  \Psi_{1,n}\left( \frac{R(z)}{1-t z^k} \right)=\left \{ \begin{array}{l} \varphi_{n-k}(R(z)), n\geq k, \\  \\ 0,  \text{ if  } n < k \end{array} \right. 
\end{array}
$$
\end{lm}
\begin{proof}

\noindent
$(i)$ Let  $h(t,z)=\sum_{i,j=0}^{\infty} h_{i,j}t^i z^j.$  Then
$$
\frac{1}{2 \pi i} \oint_{|z|=1} h\left(\frac{t}{z^n},z\right) \frac{dz}{z}=\frac{1}{2 \pi i} \sum_{i,j=0}^{\infty} \oint_{|z|=1}h_{i,j}t^i z^{j-n\,i} \frac{dz}{z}=\sum_{i}^{\infty} h_{i,ni}t^i.
$$
Thus,
$$
\frac{1}{2 \pi i} \oint_{|z|=1} h\left(\frac{t}{z^n},z\right) {\frac{dz}{z} \Bigl |_{t=z}}=\sum_{i}^{\infty} h_{i,ni}z^i= \Psi_{1,n}(h(t,z)).
$$

\noindent
$(ii)$
 Let  $R(z)=\sum_{j=0}^{\infty} f_{j} z^j.$  Then for   $k < n$  we  have 
$$
\begin{array}{l}
\displaystyle \Psi_{m,n}\left( \frac{R(z)}{1-t^m z^k} \right)=\Psi_{m,n}\Big( \sum_{j,s\geq 0} f_j  z^j (t^m z^k)^s\Big)=\Psi_{m,n}\Big(\sum_{s\geq 0} f_{s(n-k)} (t^m z^n)^s \Big){=}\sum_{s \geq 0}  f_{s(n-k)} z^s.
\end{array}
$$
On other hand, $\displaystyle \varphi_{n-k}(R(z))=\varphi_{n-k}\Bigl(\sum_{j=0}^{\infty} f_{j} z^j\Bigr){=\sum_{s \geq 0}  f_{s(n-k)} z^s.}$

 \end{proof}

The main idea of this calculations is that 
the  Poincar\'e series $ P_ {d} (z) $ can be expressed  in terms of functions $ \Psi.$ The following simple but important statement  holds
\begin{lm} 
$$
P_{d}(z)=\Psi_{1,d}\left(\frac{1+z}{(1-t)(1-t z^2)\ldots (1-t z^{2d})}\right).
$$

\end{lm}
\begin{proof} Theorem 2   implies  that   $\dim(C_{d})_n=[(t z^{d})^n]f_d(t,z^2).$ 
Then
$$
\begin{array}{l}
\displaystyle P_{d}(z) = \sum_{n=0}^{\infty}  \dim(C_{2,d})_n z^n=\sum_{n=0}^{\infty} \bigl([(t z^{d})^n]f_d(t,z^2)\bigr)z^n{=} \Psi_{1,d}(f_d(t,z^2)).
\end{array}
$$
\end{proof}
Combining this with Lemma 2, (i)  we  obtain one  more formula  for the  Puancar\'e series:

$$
\begin{array}{l}
\displaystyle P_{d}(t){=}\frac{1}{2\pi i} \oint_{|r|=1} \frac{ 1+z}{ \prod_{k=0}^{d}  (1-t  z^{d-2\,k})} \frac{dz}{z}.
\end{array}
$$
Compare the  formula  with the Molien-Weyl integral formula  for the Poincar\'e series  of the algebra of  invariants of   binary form, see \cite{DerK}, p. 183.

Now we  can present    Springer-type  formula  for the Poincar\'e  series $P_{d}(z)$
\begin{te}
$$
\begin{array}{l}
\displaystyle P_{d}(z)=\sum_{0\leq k <d/2} \varphi_{d-2\,k} \left( \frac{(-1)^k z^{k(k+1)} (1+z)}{(z^2,z^2)_k\,(z^2,z^2)_{d-k}} \right),\\
\text{ here } (a,q)_n=(1-a) (1-a\,q)\cdots (1-a\,q^{n-1}) \text{ is } q\text{-shifted factorial. }
\end{array}
$$
\end{te}
\begin{proof}
Consider the partial fraction decomposition of  the rational function $f_d(t,z^2)$$$
f_d(t,z^2)=\sum_{k=0}^{d} \frac{R_k(z)}{1-t z^{2\,k}}.
$$
It  is easy  to  see,  that 
$$
\begin{array}{l}
\displaystyle  R_k(z)=\lim_{t \to z^{-2k}}\left( f_d(t,z^2)(1-t z^k) \right)=\lim_{t \to z^{-2\,k}}\left( \frac{(1+z)}{(t,z)_{d+1}}(1-t z^{2\,k}) \right)=\\
\\
\displaystyle =\frac{1+z}{(1-z^{-2k})(1-z^{2-2k})\cdots (1-z^{2(k-1)-2k})(1-z^{2(k+1)-2k}) \cdots (1-z^{2d-2k})}=\\
\\
\displaystyle =\frac{z^{2k+(2k-2)+\ldots+2}(1+z)}{(z^{2k}-1)(z^{2k-2}-1)\cdots (z^2-1)(1-z^2) \cdots (1-z^{2d-2k}) }=

\frac{(-1)^k z^{k(k+1)}(1+z)}{(z^2,z^2)_{k} (z^2,z^2)_{d-k}}.
\end{array}
$$
Using  the above lemmas we obtain
$$
\begin{array}{l}
\displaystyle P_{d}(z) {=}\Psi_{1,d}\bigl(f_d(t,z^2)\bigr){=}\Psi_{1,d}\left( \sum_{k=0}^{n} \frac{R_k(z^2)}{1-t z^{2\,k}} \right){=} \sum_{0\leq k <d/2} \varphi_{d-2\,k} \left( \frac{(-1)^k z^{k(k+1)}(1+z)}{(z^2,z^2)_k (z^2,z^2)_{d-k}} \right).
\end{array}
$$
\end{proof}

{\bf 5.}  For direct computations we use the following technical lemma

\begin{lm} Let   $R(z)$ be some polynomial of $z.$ Then 
$$
\varphi_n\left(\frac{R(z)}{(1-z^{k_1})(1-z^{k_2})\cdots(1-z^{k_m})} \right)= 
\frac{\varphi_n\bigr(R(z)Q_n(z^{k_1})Q_n(z^{k_2})Q_n(z^{k_m})\bigr)}{(1-z^{k_1})(1-z^{k_2})\cdots(1-z^{k_m})},
$$
here  $Q_n(z)=1+z+z^2+\ldots+z^{n-1},$  and  $k_i$ are natural numbers.
\end{lm}
\begin{proof}
Observe that $ \varphi_n(f(z^n) g(z))=f(z) \varphi_n(g(z)),$ for arbitrary  $f(z),g(z) \in \mathbb{C}[[z]].$
Particularly,
$$
\varphi_n\left(\frac{g(z)}{1-z^{m}} \right)=\varphi_n\left(\frac{R(z)}{1-z^{n\,m}} \frac{1-z^{nm}}{1-z^m} \right)=\frac{1}{1-z^m}\varphi_n \left( g(z) \frac{1-z^{nm}}{1-z^m} \right)=
$$
$$
=\frac{1}{1-z^m}\varphi_n \left( g(z)(1+z^m+(z^m)^2+\cdots+(z^m)^{n-1}) \right)=\frac{1}{1-z^m}\varphi_n \left( g(z)Q_n(z^m) \right).
$$
\end{proof}

{\bf Example.}  Consider  the case $d=3.$  We  have
$$
\begin{array}{l}
\displaystyle P_3(z)=\Psi_{1,3}\left({\displaystyle \frac {1 + z}{(1-t)(1 - tz^{2})\,(1 - tz^{4})\,(1 - tz^{6}
)}}   \right)=
\\
\\
=\displaystyle \varphi_3\left(\frac{1+z}{(1-z^2)(1-z^4)(1-z^6)} \right) - \varphi_1\left({\displaystyle \frac {(1 + z)\,z^{2}}{(1
 - z^{2})^{2}\,(1 - z^{4})}} \right)=
\\
\\
=\displaystyle \frac{1}{1-z^2}\,\varphi_3\left(\frac{1+z}{(1-z^6)(1-z^{12})}\cdot \frac{1-z^6}{1-z^2} \cdot \frac{1-z^{12}}{1-z^4}\right) -\frac {(1 + z)\,z^{2}}{(1 - z^{2})^{2}\,(1 - z^{4})}=\\
\\
=\displaystyle \frac{1}{(1-z^2)^2(1-z^4)}\varphi_3\left((1+z)(1+z^2+z^4)(1+z^4+z^8) \right)-\frac {(1 + z)\,z^{2}}{(1 - z^{2})^{2}\,(1 - z^{4})}=\\
\\
=\displaystyle \frac{1+z+z^2+2z^3+z^4}{(1-z^2)^2(1-z^4)}-\frac {(1 + z)\,z^{2}}{(1 - z^{2})^{2}\,(1 - z^{4})}= \frac{1+z+z^3+z^4}{(1-z^2)^2(1-z^4)}=\\
\\
=\displaystyle  \frac{(1+z)(1+z^3)}{(1-z^2)^2(1-z^4)}= \frac{1+z^3}{(1-z)(1-z^2)(1-z^4)}.
\end{array}
$$

By using Lemma  3 the Poncar\'e series  $P_d(z)$   for  $d \leq 20$ are   found.  Below  list several results:

$$
\begin{array}{l}
\mathit{P_1(z)} := {\displaystyle \frac {1}{1-z}},
\mathit{P_2(z)} := \displaystyle \frac {1}{( - 1 + z^{2})\,(z - 1)},\\
\mathit{P_3(z)} := {\displaystyle \frac {1 + z^{3}}{(1 - z)\,(1 - z
^{2})\,(1 - z^{4})}},
\mathit{P_4(z)} := {\displaystyle \frac {1 + z^{3}}{(1 - z)\,(1 - z
^{2})^{2}\,(1 - z^{3})}} 
\end{array}
$$

\[
\mathit{P_5(z)} :={\displaystyle \frac {z^{15} + z^{13} + 3\,z^{12} + 3\,z^{11} + 5
\,z^{10} + 4\,z^{9} + 6\,z^{8} + 6\,z^{7} + 4\,z^{6} + 5\,z^{5}
 + 3\,z^{4} + 3\,z^{3} + z^{2} + 1}{(1 - z)\,(1 - z^{2})\,(1 - z
^{4})\,(1 - z^{6})\,(1 - z^{8})}}  
\]

\[
P_6(z)={\displaystyle \frac {z^{10} + z^{8} + 3\,z^{7} + 4\,z^{6} + 4\,z
^{5} + 4\,z^{4} + 3\,z^{3} + z^{2} + 1}{(1 - z)\,(1 - z^{2})^{2}
\,(1 - z^{3})\,(1 - z^{4})\,(1 - z^{5})}} 
\]

$$
\mathit{P_7(z)} := p_7(z)/((1 - z)\,(1 - z^{2})\,(1 - z^{4})\,(1 - z^{6})
 \\
(1 - z^{8})\,(1 - z^{10})\,(1 - z^{12})) 
$$
where 
$$
\begin{array}{l}
p_7(z) := z^{35} + 2\,z^{33} + 6\,z^{32} + 10\,z^{31} + 19
\,z^{30} + 28\,z^{29} + 44\,z^{28} + 61\,z^{27} + 79\,z^{26} + 
102\,z^{25} +\\
\mbox{} + 129\,z^{24} + 156\,z^{23} + 173\,z^{22} + 196\,z^{21}
 + 215\,z^{20} + 230\,z^{19} + 231\,z^{18} + 231\,z^{17}+ \\
\mbox{} + 230\,z^{16} + 215\,z^{15} + 196\,z^{14} + 173\,z^{13}
 + 156\,z^{12} + 129\,z^{11} + 102\,z^{10} + 79\,z^{9} \\
\mbox{} + 61\,z^{8} + 44\,z^{7} + 28\,z^{6} + 19\,z^{5} + 10\,z^{
4} + 6\,z^{3} + 2\,z^{2} + 1.
\end{array}
$$
$$
P_8(z)=p_8(z)/((1 - 
z)\,(1 - z^{2})^{2}\,(1 - z^{3})^{2}\,(1 - z^{4})\,(1 - z^{5})\,(
1 - z^{7})) 
$$
$$
\begin{array}{l}
p_8(z)=z^{18} + 2\,z^{16} + 6\,z^{15} + 12\,z^{14} + 19\,z^{13} + 25\,z
^{12} + 31\,z^{11} + 36\,z^{10} +38\,z^{9} + 36\,z^{8} + \\+31\,z^{
7} + 25\,z^{6} 
\mbox{} + 19\,z^{5} + 12\,z^{4} + 6\,z^{3} + 2\,z^{2} + 1
\end{array}
$$
$$
P_9(z)=p_9(z)/(1 - z)\,(1 - z^{2})\,(1 - z^{4})\,(1 - z^{6})\,(1 - z^{8})\,(1
 - z^{10})\,(1 - z^{12})\,(1 - z^{14})\,(1 - z^{16}), 
$$
$$
\begin{array}{l}
p_9(z)=z^{63} + 3\,z^{61} + 10\,z^{60} + 23\,z^{59} + 49\,z^{58} + 93\,z
^{57} + 172\,z^{56} + 289\,z^{55} + 457\,z^{54} + 701\,z^{53} \\
\mbox{} + 1036\,z^{52} + 1477\,z^{51} + 2023\,z^{50} + 2720\,z^{
49} + 3568\,z^{48} + 4573\,z^{47} + 5702\,z^{46} \\
\mbox{} + 7013\,z^{45} + 8466\,z^{44} + 10043\,z^{43} + 11672\,z
^{42} + 13400\,z^{41} + 15155\,z^{40} + 16880\,z^{39} \\
\mbox{} + 18487\,z^{38} + 20013\,z^{37} + 21392\,z^{36} + 22539\,
z^{35} + 23398\,z^{34} + 24013\,z^{33} \\
\mbox{} + 24355\,z^{32} + 24355\,z^{31} + 24013\,z^{30} + 23398\,
z^{29} + 22539\,z^{28} + 21392\,z^{27} \\
\mbox{} + 20013\,z^{26} + 18487\,z^{25} + 16880\,z^{24} + 15155\,
z^{23} + 13400\,z^{22} + 11672\,z^{21} \\
\mbox{} + 10043\,z^{20} + 8466\,z^{19} + 7013\,z^{18} + 5702\,z^{
17} + 4573\,z^{16} + 3568\,z^{15} + 2720\,z^{14} \\
\mbox{} + 2023\,z^{13} + 1477\,z^{12} + 1036\,z^{11} + 701\,z^{10
} + 457\,z^{9} + 289\,z^{8} + 172\,z^{7} + 93\,z^{6} \\
\mbox{} + 49\,z^{5} + 23\,z^{4} + 10\,z^{3} + 3\,z^{2} + 1 
\end{array}
$$
$$
P_{10}(z)=p_{10}(z)/(1 - z)\,(1 - z^{2})^{2}\,(1 - z^{3})\,(1 - z^{4})\,(1 - z^{5})\,
(1 - z^{6})\,(1 - z^{7})\,(1 - z^{8})\,(1 - z^{9}), 
$$
$$
\begin{array}{l}
p_{10}(z)=z^{36} + 3\,z^{34} + 11\,z^{33} + 27\,z^{32} + 58\,z^{31} + 112\,
z^{30} + 193\,z^{29} + 318\,z^{28} + 485\,z^{27} + 699\,z^{26}
 \\
\mbox{} + 951\,z^{25} + 1245\,z^{24} + 1541\,z^{23} + 1842\,z^{22
} + 2108\,z^{21} + 2321\,z^{20} + 2451\,z^{19} \\
\mbox{} + 2506\,z^{18} + 2451\,z^{17} + 2321\,z^{16} + 2108\,z^{
15} + 1842\,z^{14} + 1541\,z^{13} + 1245\,z^{12} \\
\mbox{} + 951\,z^{11} + 699\,z^{10} + 485\,z^{9} + 318\,z^{8} + 
193\,z^{7} + 112\,z^{6} + 58\,z^{5} + 27\,z^{4} + 11\,z^{3} \\
\mbox{} + 3\,z^{2} + 1.
\end{array}
$$


\begin{thebibliography}{15}


\bibitem{SF}
Sylvester, J. J.,  Franklin, F., Tables of the generating functions and groundforms for the binary quantic of the first ten orders, Am. J. II., 223--251, 1879
\bibitem{Sylv-12}
Sylvester, J. J.
Tables of the generating functions and groundforms of the binary duodecimic, with some general remarks, and tables of the irreductible syzigies of certain quantics. Am. J. IV. 41-62, 1881.



\bibitem{SP}
Springer T.A., On the invariant theory of SU(2), Indag. Math. 42 (1980), 339-345.

\bibitem{BC}    Brouwer A.,  Cohen A., The Poincare series of the polynomial invariants under
$SU_2$  in its irreducible representation of degree $\leq  17$, preprint of the Mathematisch Centrum,
Amsterdam, 1979.

\bibitem{LP}
 Littelmann P.,  Procesi C., On the Poincar\'e series of the invariants of binary forms, J. Algebra 133, No.2, 490-499 (1990).

\bibitem{Rob}
Roberts M., The covariants of a binary quantic of the n-th degree,
Quarterly J. Math., {\bf 4},1861,168--178.

\bibitem{FH} W.Fulton,  J. Harris, Reptesentation theory: a first course, 1991.

\bibitem{DerK}
 Derksen H.,  Kemper G., Computational Invariant Theory, Springer-Verlag, New York,
2002.
\end{thebibliography}
\end{document}